\documentclass{amsart}
\usepackage[mathscr]{eucal}
\usepackage{amssymb, amsmath,array, amscd}
\usepackage{enumerate}
\usepackage{graphicx}
\usepackage{url}
\usepackage[colorlinks,plainpages,backref]{hyperref}

\input xy
\xyoption{all}

\newtheorem{THM}{Theorem} [section]

\newtheorem{prop}{Proposition}[section]
\newtheorem{lemma}{Lemma}[section]

\newtheorem{remark}{Remark}[section]
\newtheorem{cor}{Corollary}[section]

\newtheorem{example}{Example}[section]

\begin{document}

\title{Completely reducible hypersurfaces in a pencil}

\author{J. V. Pereira}
\address{IMPA \\ Est.
D. Castorina, 110\\
22460-320, Rio de Janeiro, RJ, Brasil} \email{ jvp@impa.br}

\author{S. Yuzvinsky}
\address{Departament of Mathematics \\ University of Oregon \\ Eugene \\ OR 94703 USA}
\email { yuz@math.uoregon.edu }

\dedicatory{}
\thanks{The first author is supported by Cnpq and Instituto Unibanco.}

\
\begin{abstract}
We study completely reducible fibers of pencils of hypersurfaces
on $\mathbb P^n$ and associated codimension one foliations of $\mathbb P^n$.
 Using methods from
theory of foliations we obtain certain upper bounds for the
number of these fibers as functions only of $n$.
 Equivalently this
gives upper bounds for the
dimensions of resonance varieties of hyperplane arrangements.
 We obtain similar bounds
for the dimensions of the characteristic varieties of
the arrangement complements.
\end{abstract}

\maketitle

\date{\today}

\section{Introduction}
\label{S:1}

{In this paper}  we focus on completely reducible fibers of a pencil
of hypersurfaces on $\mathbb P^n$ with irreducible generic fiber.
Our main result (Theorem \ref{T:0}) gives an upper bound for the
number $k$ of these fibers that depends only on $n$. For instance
for every $n>1$ we obtain $k\le 5$. Or for every $n\ge 4$ a pencil
with $k$ at least 3 is a linear pull-back of a pencil on $\mathbb
P^4$.

Despite the evident classical taste of the result
we have not found it in the literature although various restrictions
on reducible fibers and even more special completely reducible fibers
were studied in \cite{Vi, Had, Hal}.

In order to prove the result we use a combination of techniques from
{the} theory of codimension one singular foliations on $\mathbb P^n$
and {the} theory of hyperplane arrangements. In particular we define
a foliation associated with a pencil and consider its Gauss map
  $$\mathbb P^n \dashrightarrow (\mathbb P^n)^{\vee}. $$
 Our most technical
result (Theorem \ref{T:1}) says that for a pencil with $k\ge 3$
the Gauss map of the associated foliation is dominant.

It turns out that the union of the linear divisors of all completely
reducible fibers of a pencil with $k$ at least three can be characterized
intrinsically, by means of arrangement theory. In fact this characterization
in terms are resonance varieties of hyperplane arrangements
is contained in the recent paper \cite{FY} and this was the starting point of our study.
We recall this characterization in a concise form in Theorem \ref{T:5}.
This and Theorem  \ref{T:0} give upper bounds on the dimensions of the resonance varieties
of arrangements whence on the dimension of cohomology of the Orlik-Solomon algebras.

Much of the interest in the resonance varieties comes from the fact
that they are closely related to the support loci for the cohomology
of local systems on arrangement complements - the so called
characteristic varieties. These relations have been studied by many
authors, see for instance \cite{CS1999},\cite{LY2000}. In a recent
preprint \cite{Dimca} Dimca cleared some subtle points in these
relations {as well as} in the relations of characteristic components
with pencils. Using this we are able to find upper bounds for
dimensions of these components too (Theorem \ref{T:7}). In order to
do this in full generality (at least {for} the {positive}
dimensional components) we need to generalize Theorem \ref{T:1} to
more general pencils having at least two completely reducible fibers
and another fiber with irreducible components either of degree one
or non-reduced. Since this condition on pencils does not look
natural we does not include this generalization (Theorem \ref{T:6})
in Theorem \ref{T:1}.

We do not know if in all the cases the upper bounds given in Theorems \ref{T:1}
and \ref{T:6} are strict. In section \ref{S:6} we collect some old and new examples
of pencils on $\mathbb P^n$ ($n>1$)
with at least three completely reducible fibers.

\section{Main result and reduction}
\label{S:2}

\subsection{Completely reducible fibers}
Let $F$ and $G$ be polynomials of the same degree $d>0$ from
$\mathbb C[x_0,x_1,\ldots,x_n]$ defined up to nonzero multiplicative
constants. {They} define the pencil $\mathcal P=\{aF+bG\, |\,
[a:b]\in\mathbb P^1\}$ whose fibers $aF+bG$ can be identified with
hypersurfaces in $\mathbb P^n$. We are looking for upper bounds for
the number of the completely reducible fibers, i.e., products of
linear forms. Without loss of generality we can and will always
assume that $F$ and $G$ are relatively prime (equivalently, the
generic fiber of $\mathcal P$ is irreducible) and $F$ and $G$ are
completely reducible themselves, whence there are at least two
completely reducible fibers of $\mathcal P$. In the trivial case
where $d=1$, all fibers of the pencil $\mathcal P$ are hyperplanes
whence completely reducible.

The following theorem is the main result of the paper.

\begin{THM}\label{T:0}
If $\mathcal P$ is a pencil of hypersurfaces on $\mathbb P^n$ with
irreducible generic fiber and $k$ is the number of completely
reducible fibers of $\mathcal P$ then the following assertions hold
\begin{enumerate}
\item If $k > 5 $ then $\mathcal P$ is a pencil of hyperplanes (equivalently, it is
the linear pull-back of a pencil on $\mathbb P^1$);
\item If $k > 3$ then $\mathcal P$ is the linear pull-back of a
pencil on $\mathbb P^2$;
\item If $k > 2$ then $\mathcal P$ is the linear pull-back of a
pencil on $\mathbb P^4$.
\end{enumerate}
\end{THM}

\smallskip

\subsection{Reduction}
\label{reduction} The proof of Theorem \ref{T:0} {consists} of two
parts. In this subsection we present an elementary reduction to the
non-vanishing of {a} certain determinant. The rest of the proof
requires more techniques and will be given in the following
sections.

\medskip

We are using the notation from the theorem. If $k\le 2$
there is nothing to prove. So we assume that $k\ge 3$.

Let $\widetilde Q=\prod\alpha_i^{m_i}$ denotes the product of the
completely reducible fibers of $\mathcal P$ and $Q=\prod\alpha_i$
denotes the reduced polynomial with the same zero set. Then consider
the $1$-form
\[
 \omega = \frac{Q}{\widetilde Q} \omega_0
\]
where $\omega_0 = FdG -GdF$.

We will need the following properties of $\omega$.

(i) \emph{$\omega$ does not depend (up to a nonzero multiplicative constant)
on the choice of $F$ and $G$ from the set of completely reducible fibers of $\mathcal P$}.

A proof is left to the reader.

(ii) \emph{$\omega$ is a polynomial form}.

It suffices to check that $\alpha_i^{m_i-1}$ (for every $i$) divides the coefficients of
$\omega_0$. If $\alpha_i$ divides $F$ (or $G$) this is clear; otherwise it follows from (i).

For the future use we write $\omega=\sum_{i=1}^{n=1}a_idx_i$ for some $a_i\in\mathbb C[x_0,
\ldots,x_n]$ and denote by $D$ the determinant of the Jacobi matrix
$\left(\frac{da_i}{dx_j}\right)$.

(iii) \emph{For every linear factor $\alpha_i$ of $Q$
we have $\alpha_i$ divides the coefficients
of $d\alpha_i\wedge\omega$}.

\begin{proof} It suffices to check that $\alpha_i^{m_i}$ divides $d\alpha_i\wedge\omega_0$.
Again if $\alpha_i$ divides $F$ one needs to check that $\alpha_i^{m_i}$
divides $d\alpha_i\wedge dF$ which is clear. Otherwise one can apply (i).
\end{proof}

\medskip

(iv) \emph{$Q^{n-1}$ divides $D$}.

\begin{proof}
It suffices to prove that $\alpha_i^{n-1}$ divides $D$ for every $i$.
Fix $i$ and change the coordinates so that $\alpha_i=x_0$.
Applying property (iii) we see that $x_0$ divides $a_j$ for $j=1,2,\ldots,n$.
In particular, except for the entries of the first row and the first column,
all the entries of the matrix $\left(\frac{\partial a_j}{\partial x_k}\right)$
are divisible by $x_0$. Using the cofactor expansion of $D$ with respect to the first
row completes the proof.
\end{proof}

Now we can prove that the non-vanishing of $D$ would imply Theorem \ref{T:0}.

\begin{prop}
\label{red}
Suppose $D$ is not identically 0. Then Theorem \ref{T:0} follows.
\end{prop}
\begin{proof}

It follows from property (iv) of $\omega$ that  $\widetilde Q^{n-1}$
divides
\[
  \left(\frac{\widetilde Q}{Q}\right)^{n-1}D
  \, .
\]
 Since $\deg(\widetilde Q) = k \deg(F)=kd$ we
obtain that
\begin{eqnarray*}
(n-1) kd  &\le& (n+1) \left(2d  -2 -\deg\left(\frac{\widetilde
Q}{Q}\right)\right) +
(n-1)\deg\left(\frac{\widetilde Q}{Q}\right) \\
&=& (n+1) \left(2d  -2 \right)  - 2\deg\left(\frac{\widetilde
Q}{Q}\right).
\end{eqnarray*}
Therefore
\[
k \le \frac{n+1}{n-1} \left(2  - \frac{2}{d} \right) -
\frac{2}{(n-1)d}\deg\left(\frac{\widetilde Q}{Q}\right) <
2\frac{n+1}{n-1} \, .
\]
In particular
\begin{eqnarray*}
n\ge 5 & {\rm implies} & k \le 2 \,  , \\
n\ge 3 & {\rm implies} & k \le 3 \, , \\
n\ge 2 & {\rm implies} & k \le 5 \, .
\end{eqnarray*}
Theorem \ref{T:0} follows.
\end{proof}

To prove that the polynomial $D$ is not 0 we will use foliations on
$\mathbb P^n$ and some properties of the hyperplane
multi-arrangement $\widetilde Q=0$.

\medskip

\section{The multinet property and its generalization}\label{S:multinet}
\label{S:3}

First in this section we recall certain results from \cite{FY} about
the multi-arrangement $(\mathcal A,m)$ defined by $\widetilde Q=0$.
We remark that we recall here  only {the} facts about arrangements
that are needed for the proof of Theorem \ref{T:0}. For more
{details} see section \ref{S:5}.

The multi-arrangement $(\mathcal A,m)$
 consists of the set $\mathcal A$ of hyperplanes $H_i$ in $\mathbb P^n$
 determined by $\alpha_i=0$
where $\alpha_i$ is running through all linear divisors of $Q$.
Besides to each $H_i$ the positive integer $m(H_i)=m_i$ is assigned
where $m_i$ is the exponent of $\alpha_i$ in $\widetilde Q$.

First we {consider} the case where $n=2$. In this case, it was
proved in \cite{FY} that the partition
\[
\mathcal A = \mathcal A_1 \cup \mathcal A_2 \cup \cdots \cup
\mathcal A_{k}
\]
of $\mathcal A$ into completely reducible fibers (called classes in \cite{FY})
of the pencil $\mathcal P$
can be equivalently characterized by the combinatorics of lines and points.
In particular
\begin{itemize}
\item[(a)] $\sum_{H \in \mathcal A_i} m(H)$ is independent of
$i=1, \ldots, k$;
\item[(b)] If $p$ is the point in the base locus of $\mathcal P$
then the sum
\[
  n(p)=\sum_{H \in \mathcal A_i, p \in H} m(H)
\]
is independent of $i= 1, \ldots, k$.
\end{itemize}

The collection $(\mathcal A,\mathcal X)$ of lines and points
satisfying conditions (a)-(d) is called a \emph{k-multinet}. If
$n(p)=1$ for all $p\in\mathcal X$ (whence $ m(H)=1$ for all
$H\in\mathcal A$) then it is a \emph{k-net}.

Intersecting an arrangement in $\mathbb P^n$ with a general position
plane $\mathbb P^2$ one readily sees that
the similar properties hold for arrangements of hyperplanes in
$\mathbb P^n$ where the intersections of codimension 2 should be substituted
for points. This definition of multinets in $\mathbb P^n$ for $n>2$
does not say anything about the intersections of higher codimensions of
hyperplanes in a multinet. We can prove and then use a property of these intersections
for $k\ge 3$. In fact we prove this property under a weaker assumption
that we will use in section \ref{S:7}.

\begin{prop}\label{L:like}
Let $\mathcal P$ be a pencil of hypersurfaces on $\mathbb P^n$ with
irreducible generic fiber generated by two completely reducible
fibers $F$ and $G$. If there exists a third fiber that is a product
of linear forms  and non-reduced polynomials then
\[
\sum_{\alpha_H|F, p \in H} m(H) = \sum_{\alpha_H|G, p \in H} m(H)
\]
for every $p$ in the base locus of $\mathcal P$.
\end{prop}
\begin{proof}
Let $p$ be a point in the base locus of the pencil.
Choose affine coordinates $(x_1,\ldots,x_n)$ where $p$ is the origin and write
$F=F_1\cdot F_2$ and $G=G_1 \cdot G_2$ where $F_2,G_2 \notin \mathfrak m$ and all the irreducible
components of $F_1$ and $G_1$ are in $\mathfrak m$,
$\mathfrak m$ being the maximal ideal $(x_1,\ldots, x_n)$.
 The statement of the lemma is equivalent to {the} homogeneous polynomials
$F_1$ and $G_1$ having the same {degree}.

Our hypothesis implies that there exists a hyperplane or a
non-reduced  hypersurface {in the pencil} passing through $0$. In
the former case, there is a linear form $\alpha\in{\mathbb
C}[x_1,\ldots,x_n]$ that divides, say $K=F-G$. Thus $\alpha$ divides
$K_0=F_2(0)F_1-G_2(0)G_1$. If $\deg F_1\not=\deg G_1$ then $\alpha$
divides, say $F_1$ which is a contradiction.

In
the latter case, there is an irreducible polynomial $f\in {\mathfrak m}$
 such that $f^m$ divides $K$ for some
$m>1$ whence $f$ divides the coefficients of $\omega=FdG-GdF=KdF-FdK$.
Now put $R=x_1\frac{\partial}{\partial x_1} + \cdots + x_n\frac{\partial}{\partial x_n}$,
i.e., $R$ is the radial (or Euler) vector field in $\mathbb C^n$, and denote by $i_R$
the interior product of a form with it.
Then the Leibniz formula implies that
\[
i_R\omega =  F_2G_2 i_R({\underbrace{F_1dG_1-G_1dF_1}_{\omega_1}}) +
F_1G_1 i_R \left(F_2dG_2 - G_2dF_2\right) .
\]
To prove the lemma it suffices to show that $i_R\omega_1=0$.

Suppose this is not true. Then we have
$i_R\omega_1=cF_1G_1$ for a $c\in\mathbb{C^*}$ and
$i_R\omega=F_1G_1g$ where $g$ is a polynomial in $x_1,\ldots,x_n$
such that $g(0)\not=0$.
Since $f$
divides the polynomial $i_R\omega$ it divides $F_1G_1$
which is again a contradiction.
\end{proof}

Sometimes it is convenient to assume that the pencil $\mathcal P$ is not a linear pull-back
from a smaller dimensional projective space. We will say in this case
that $\mathcal P$ is \emph{essential} in $\mathbb P^n$. This is equivalent
to $\mathcal A$ being \emph{essential}, i.e., $\bigcap_{H\in\mathcal A}H=\emptyset$.
This can be
expressed also by saying that the rank of $\mathcal A$
 is $n$ where rank is the codimension of $\bigcap_{H\in\mathcal A}H$
(see section \ref{S:5}).
It is immediate from the multinet property that $\mathcal A$ is essential if and only if
the arrangement defined by all the linear divisors of $F$ and $G$ is essential.

We will need the following property of essential arrangements that immediately follows from
definitions.

\begin{prop}
\label{property}

Let $\mathcal A$ be the collection of the linear divisors of all the completely
reducible fibers of an essential pencil $\mathcal P$ on $\mathbb P^n$.
Then there exist two distinct points $p_1,p_2\in\mathbb P^n$ such that
for $j=1,2$  the subarrangements
   $$\mathcal B_j = \bigcup_{p_j \in H, H \in \mathcal A} H$$
have rank $n-1$.
\end{prop}

\section{Foliations and the Gauss map}
\label{S:4}

\subsection{Foliations.}

In this paper we will adopt {an} utilitarian definition for
codimension one singular foliations  on $\mathbb P^n$, from now on
just foliation on $\mathbb P^n$. A \emph{foliation} $\mathcal F$ on
$\mathbb P^n$  will be an equivalence class of homogeneous rational
differential $1$-forms on $\mathbb C^{n+1}$ under the equivalence
relation
\[
\omega \sim \omega ' \, \, \text{ if and only if there exists } h \in \mathbb C[x_0,
\ldots, x_{n}]\setminus 0 \, \text{ for which  } \, \omega = h
\omega',
\]
such that $i_R \omega =0$  and $\omega \wedge d\omega =0 $ for every
 representative $\omega$. Here $R$ and $i_R$ are similar to the ones used in the
proof of Proposition \ref{L:like} but in $\mathbb C^{n+1}$.
 Of course, to
ensure the validity of the two conditions for every representative it is
sufficient to check it just for one of them.

Among the representatives of $\mathcal F$ there are privileged ones
--- the homogeneous polynomial $1$-forms with singular, i.e.
vanishing, set of codimension at least two. Any two such forms that
are equivalent differ by a nonzero multiplicative constant. If such
a form has coefficients of degree $d+1$ then we say that $\mathcal
F$ is a degree $d$ foliation. The shift in the degree is motivated
by the geometric interpretation of the degree. It is the number of
tangencies between $\mathcal F$ and a generic line in $\mathbb P^n$.

Outside the singular set the well-known Frobenius Theorem ensures
the existence of  local submersions  with connected level sets
whose tangent space at a point is the kernel of a defining $1$-form
at this point. These level sets are the local leaves of $\mathcal F$.
The leaves are obtained by patching together level sets of distinct
submersions that have nonempty intersection. Although the data is
algebraic the leaves, in general, have a transcendental nature.

Now we show how foliations appear from completely reducible fibers
of a pencil of {hypersurfaces} on $\mathbb P^n$.

\begin{lemma}
\label{foliation}
Let $\mathcal P$ be the pencil on $\mathbb P^n$ generated by polynomials $F$ and $G$
and $\omega$ the 1-form from \ref{S:2}. Then $\omega$ defines a
foliation on $\mathbb P^n$.
\end{lemma}
\begin{proof}
One needs to check the two conditions from definition of foliations
for $\omega$ or equivalently, for the closed form
$\eta=\frac{dG}{G}-\frac{dF}{F}$. The condition $i_R\eta=0$ follows
immediately since $F$ and $G$ are homogeneous polynomials of equal
degrees, {while the integrability condition is automatically
satisfied thanks to the closedness of $\eta$}.
\end{proof}

The foliation defined by $\omega$ will be called the foliation {\it
associated to $\mathcal P$}.

\subsection{Gauss map}
Let $\omega$ be a  homogeneous polynomial differential $1$-form on $\mathbb C^{n+1}$ such
that $i_R\omega =0$ and $\omega\wedge d \omega =0$. Let $\mathcal F$
be the foliation defined by $\omega$ on $\mathbb P^n$. {\emph The Gauss map} of
$\mathcal F$ is the rational map
\begin{eqnarray*}
\mathcal
 G_{\omega}= G_{\mathcal F} : \mathbb P^n &\dashrightarrow& (\mathbb
  P^n)^{\vee} \\
  p &\mapsto& T_p\mathcal F \, .
\end{eqnarray*}
that takes every point $p \in \mathbb P^n \setminus
\mathrm{sing}(\mathcal F)$ to the hyperplane tangent to $\mathcal F$
at p. Under a suitable identification of $\mathbb P^n$ with
$(\mathbb P^n)^{\vee}$ the Gauss map $\mathcal G_{\omega}$  is
nothing more than the rational map defined in homogeneous
coordinates by the coefficients of $\omega$.

\smallskip

We say that a foliation $\mathcal F$ has degenerate Gauss map when
$\mathcal G_{\mathcal F}$ is not dominant, i.e., its image is not
dense in $(\mathbb P^n)^{\vee}$. On $\mathbb P^2$ a
foliation with degenerate Gauss map has to be a
pencil of lines. Indeed the restriction of the Gauss map to a leaf
of the foliation coincides with the Gauss map of the leaf and the
only (germs of) curves on $\mathbb P^2$ with degenerate Gauss map
are (germs of) lines. Thus  all the leaves of a foliation with
degenerate Gauss map are open subsets of lines. In order for this foliation not
to have every point singular, these
lines should intersect at one point.

On $\mathbb P^3$ the situation is more subtle and a complete
classification can be found in \cite{CeLN}. Some results toward the
classification of foliations with degenerate Gauss map on $\mathbb
P^4$ have been recently obtained by T. {Fassarella} \cite{thiago}.

In our special case where a foliation is associated to the pencil $\mathcal P$
we can prove that the Gauss map cannot degenerate.
\medskip

\begin{THM}\label{T:1}
Let $\mathcal P$ be an essential pencil of hypersurfaces on $\mathbb P^n$
with at least three completely reducible fibers.
Then the Gauss map of the associated
foliation $\mathcal F$ is non-degenerate.
\end{THM}
\begin{proof}
We again denote by $\mathcal A$ the arrangement defined by the
linear divisors of all $k$ completely reducible fibers of $\mathcal
P$ and let $\mathcal A=\mathcal A_1\cup\cdots \cup\mathcal A_k$ be
its partition into fibers. We assume that $\mathcal A_1$ and
$\mathcal A_2$ correspond to $F$ and $G$ respectively. For each
$H\in\mathcal A$ we denote by $m(H)$ the multiplicity of the
respective linear form $\alpha_H$ in the its fiber.

Now we use induction on $n$. Suppose $n=2$. The only way to have a pencil of lines
as the foliation
associated to $\mathcal P$ is for $\mathcal A$ to be a pencil of lines itself. But
a pencil of lines is not essential in $\mathbb P^2$. Thus the
Gauss map is not degenerate in this case.

Suppose that the result holds for essential pencils in
$\mathbb P^n$ and let $\mathcal P$ be  an essential pencil
in $\mathbb P^{n+1}$.
 The foliation associated to $\mathcal P$ can be thus
defined by the $1$-form
\[
\omega_0= FdG - GdF \ ,
\]
where
\[
F = \prod_{H \in \mathcal A_1} {\alpha_H}^{m(H)} \quad \text{and}
\quad G = \prod_{H \in \mathcal A_2} {\alpha_H}^{m(H)}  \, .
\]

Fix two points
$p_1,p_2 \in \mathbb P^{n+1}$ with the property from Proposition \ref{property} and
let $\pi_1 : \widetilde{\mathbb P^{n+1}} \to \mathbb P^{n+1}$ be the
blow-up of $\mathbb P^{n+1}$ at $p_1$. Using notation from Proposition
\ref{property} the restriction to exceptional divisor  $E_1\cong
\mathbb P^{n}$ of the strict transforms of the hyperplanes in
$\mathcal B_1$ induces a non-degenerate arrangement of hyperplanes
in $E_1$ which we will still denote by $\mathcal B_1$.

If $\iota: E_1  \to  \widetilde{ \mathbb P^{n+1}}$ is the natural
inclusion then we claim that  the closure of the image of the
rational map $\sigma_1= \mathcal G_{\mathcal F} \circ \pi_1 \circ
\iota$ (see the diagram) is the hyperplane $H^{(1)}$ in $({\mathbb P^{n+1}})^{\vee}$ dual to
$p_1$.
\[
\xymatrix{
E_1\cong \mathbb P^n \ar[r]^{\iota}\ar@{.>}@/^1.5cm/[rrrd]^{\sigma_1}  & \widetilde{\mathbb P^{n+1}} \ar[d]_{\pi_1}\ar@{.>}[drr]^{ \mathcal G_{\mathcal F} \circ \pi_1}\\
 & \mathbb P^{n+1} \ar@{.>}[rr]^-{\mathcal G_{\mathcal F}} & & {(\mathbb
P^{n+1})}^\vee}
\]

Indeed the arrangement $\mathcal B_1 \subset E_1$ admits a partition
\[
 \mathcal B_1 = \bigcup_{i=1}^{k} \mathcal B_1 \cap \mathcal A_i
\]
and a function $m_1 = m_{\vert \mathcal B_1}$ satisfies the
multinet properties with the initial multiplicities and the induced partition into
classes. This follows from Propositions \ref{L:like} and \ref{property}.
Now the equivalence from section \ref{S:3} implies that
$\prod_{H\in\mathcal A_i\cap\mathcal B_1}\alpha_H^{m(H)}$
is a fiber of the pencil $\mathcal P_1$ on $E_1$ generated by
\[
F_1 = \prod_{H \in \mathcal A_1\cap \mathcal B_1 } \alpha_H^{m(H)}
\quad \text{and} \quad  G_1 = \prod_{H \in \mathcal A_2\cap \mathcal
B_1} \alpha_H^{m(H)} \, .
\]

Now consider the associated foliations.
On one hand the foliation $\mathcal F_1$
on $\mathbb P^n$ associated to $\mathcal P_1$
can be defined by the $1$-form $F_1dG_1 - G_1dF_1$.

On the other hand the first nonzero jet of $FdG-GdF$ at $p_1$ is
$F_1dG_1 - G_1dF_1$ viewing now  $\alpha_H$ as linear forms on
$\mathbb C^{n+2}$. More precisely let us assume that
$p_1=[0:\ldots:0:1]$ and write $F=F_1F_2, \ G=G_1G_2$. Then the
map $\mathcal G_{\mathcal F}$ (after division by $F_2(p_1)G_2(p_1)$) can be
written in the homogeneous coordinates
$[x_0:x_1:\ldots:x_n:x_{n+1}]$ of $ \mathbb P^{n+1}$ as
\begin{eqnarray*}
\mathcal G_{\mathcal F} &=& \left[ F_k\frac{\partial G_1}{\partial
x_0}-G_1 \frac{\partial F_1}{\partial x_0} + b_0 : \cdots:
F_1\frac{\partial G_1}{\partial x_n}-G_1 \frac{\partial
F_1}{\partial x_n} + b_n : b_{n+1} \right] \\
&=& \left[F_1dG_1 - G_1dF_1 + \sum_{i=0}^{n+1} b_i dx_i \right],
\end{eqnarray*}
where $b_i\in\mathfrak m^{2d}$ ($d=\deg F$) for $i=0,\ldots,n,n+1$,
$\mathfrak m$ being the maximal ideal of $\mathbb C[x_0,\ldots,
x_n]$ supported at $0 \in \mathbb C^{n+1}$.

If we consider now $(x_0:\ldots:x_n) \in \mathbb P^n$ as a
homogenous system of coordinates on the exceptional divisor $E_1$
then in these {coordinates}
\[
    \sigma_1 = \left[ F_1\frac{\partial G_1}{\partial
x_0}-G_1 \frac{\partial F_1}{\partial x_0} : \cdots:
F_1\frac{\partial G_1}{\partial x_n}-G_1 \frac{\partial
F_1}{\partial x_n}  \right] = [F_1dG_1 - G_1dF_1]
\]
since the coefficients of  $F_1dG_1 - G_1dF_1$ lie in
$\mathfrak m^{2d-1}$.  Thus $\sigma_1$ can be identified with the
Gauss map of the foliation $\mathcal F_1$
 composed with $\phi$, where $\phi$ is
the isomorphism of $({\mathbb P^n})^{\vee}$
with $H^{(1)}$ defined by the coordinates.

Using also the point $p_2$ and applying
the inductive hypothesis we have now that the closure of the image of
$\mathcal G_{\mathcal  F}$ contains at least two distinct
hyperplanes. Since $\mathbb P^{n+1}$ is irreducible  the closure of
the image of $\mathcal G_{\mathcal F}$ must be  $({\mathbb
P^{n+1}})^{\vee}$. This completes the proof.
\end{proof}

Now Theorem \ref{T:1} and Proposition \ref{red} constitute a proof of Theorem \ref{T:0}.

\medskip

\subsection{Invariant hyperplanes}

As another application of the Gauss map we can exhibit an upper bound
on the number of hyperplanes invariant with respect to a foliation in $\mathbb P^n$.

Let $\mathcal F$ be a foliation on $\mathbb P^n$ defined by a polynomial
1-form $\omega=\sum_{i=0}^na_idx_i$.  Recall that the invariance of a hyperplane $H$
defined by a linear form $\alpha_H$ with respect to $\mathcal F$ means that
$\alpha_H$ divides all coefficients of the form $\omega\wedge d\alpha_H$, i.e.,
the property (iii) from subsection \ref{reduction} holds for $\omega$ and $\alpha_H$.
Then the property (iv) holds also, i.e., $\alpha_H^{n-1}$ divides the
determinant $D$ of the matrix $\left(\frac{\partial a_i}{\partial x_j}\right)$.
Taking as $\omega$ a polinomial form without a codimension one zero we have
$\deg D=(n+1)\deg{\mathcal F}$. This gives the following proposition.

\begin{prop}\label{P:para}
If $\mathcal F$ is a foliation on $\mathbb P^n$ with
the non-degenerate Gauss map then the number of invariant hyperplanes is
at most
\[
\left( \frac{ n+1} { n-1}  \right) \cdot \deg (\mathcal F) \, .
\]
\end{prop}

\medskip

\begin{example}\rm
Here is an example showing that the bound in Proposition \ref{P:para}  is sharp.
Let $\mathcal F$ be a foliation on $\mathbb P^n$, $n \ge 2$,
induced by a logarithmic $1$-form
\[
 \sum_{i=0}^n \lambda_i \frac{dx_i}{x_i}
\]
where $\sum \lambda_i=0$ and no $\lambda_i$ is equal to zero.
Then $\mathcal F$ has degree $(n-1)$ and
its Gauss map is non-degenerate. Also $\mathcal F$ leaves invariant all the
$n+1$  hyperplanes of the arrangement.
\end{example}

For $n=2,3$ there are  examples of foliations
$\mathcal F$ on $\mathbb P^n$ with $\deg\mathcal F>n-1$ and
 exactly $\frac{ n+1} { n-1} \deg
(\mathcal F)$ invariant hyperplanes. On $\mathbb P^2$ we are aware
of three sporadic examplesi: Hesse pencil, Hilbert modular foliation
\cite{lg}, and a degree $7$ foliation leaving invariant the extended
Hessian arrangement of all the reflection hyperplanes of the reflection
group of order 1296 - see \cite{OT}, p. 227. Also there is one infinite family \cite{jvp},
\cite[Example 4.6]{FY},  consisting
of degree $m$, $m\ge 2$, foliations leaving invariant the
arrangement
\[
x y z ( x^{m-1} - y^{m-1} ) ( x^{m-1} - z^{m-1} ) ( z^{m-1} -
y^{m-1} ).
\]

On $\mathbb P^3$ we are aware of just one example with $\deg\mathcal F> n-1$
attaining the bound, see \ref{S:P3} below.

\section{Characterization of the union of completely reducible fibers}
\label{S:5}

\subsection{Hyperplane arrangements.}
First we need to recall some facts about hyperplane arrangements.
Let $\mathcal A = \{ H_1, \ldots, H_m \}$ be an arrangement of
linear hyperplanes in $\mathbb C^{n}$. Recall that the rank of $\mathcal A$ is
the codimension of $\bigcap_iH_i$. In particular if the rank is $n$
the arrangement is essential. Let
$M = \mathbb C^{n} \setminus \bigcup_{H\in\mathcal A}H$ be the
complement of $\mathcal A$. The cohomology ring $H^*(M)$ is well-known (for instance,
see \cite{OT}).
In particular it does not have torsion and working with complex coefficients
does not loose any generality.

Put $A=H^*(M,\mathbb C)$ and as before fix for each $i=1,2,\ldots,m$ a linear form
$\alpha_i$ with the kernel $H_i$. By the celebrated Arnold-Brieskorn theorem
(e.g., see \cite{OT}) the algebra $A$ is isomorphic under the deRham map
to the subalgebra of the algebra of the closed differential forms
generated by $\{\frac{d\alpha_i}{\alpha_i}, 1=1,2,\ldots,m\}$.
Notice that $A^1$ is a linear space of dimension $m$.

Since $A$ is graded commutative each $a \in A^1$ induces a cochain complex
\[
(A,a): \quad 0 \to A^0 \to A^1 \to A^2 \to \cdots \to A^{k} \to
A^{k+1} \to 0
\]
where the differential is defined as the multiplication by $a$. The degree
$l$ \emph{resonance variety} $\mathcal R^l(\mathcal A)$ is
\[
\mathcal R^l(\mathcal A ) = \{ [a] \in \mathbb P( A^1)\cong \mathbb
P^{m-1} \, \vert \, {\mathop H}^l(A,a) \neq 0 \} \, .
\]

In this paper we will need only
the first resonance variety $\mathcal R^1(\mathcal A )$.
If $\mathcal A' \subset \mathcal A$ is a subarrangement
then by definition $\mathcal R^1(\mathcal A' ) \subset \mathcal R^1(\mathcal A )$.
The \emph{support} of  an irreducible component $\Sigma$ of $\mathcal
R^1(\mathcal A )$ is the smallest subarrangement $\mathcal A'
\subset \mathcal A$ such that  $\Sigma \subset \mathcal R^1(\mathcal
A' )$. For us the rank of the support of an irreducible component is
important. The irreducible components will be called \emph{global} if
the rank of its support equals the rank of
$\mathcal A$. In the rest of the paper we call irreducible components
of $\mathcal R^1$ \emph{resonance components}. It is well-known
(see \cite{CS1999}) that the resonance components are linear subspaces of $A^1$.

Since we study pencils in projective spaces we need to projectivize linear
arrangements, i.e., to deal with arrangements of hyperplanes in the
projective space $\mathbb P^n$. We still call it essential if its
linear cone is such. More explicitly this means that the intersection of all
hyperplanes is empty. According to Proposition \ref{property}
there are at least two points among intersections
of hyperplanes and the subarrangement of
the hyperplanes passing through any of those points is isomorphic to
a linear arrangement of rank $n$.
We say that the rank of the essential projective arrangement is $n$.
More generally the rank of an arbitrary projective arrangement is the rank of its
linear cone minus one.
For the resonance varieties of a projective arrangement
 we still use those of the linear cone of it.

\subsection{Pencils and resonance components}
In this subsection we give a characterization of
the unions of completely reducible fibers of pencils of
hypersufaces on $\mathbb P^n$ in terms of resonance components.
We also interpret the partition
of such an arrangement into fibers and give a corollary of our main result
for the dimensions of the resonance components of arrangements.

The first part is not really new; the following theorem is just a reinterpretation
of results from \cite{FY}, in particular Corollary 3.12.
Although these results have been proved there only for
pencils in $\mathbb P^2$ they can be immediately generalized to arbitrary
dimensions using intersection with a general position plane.

\begin{THM}\label{T:5}
Let $\mathcal A$ be a projective arrangement in $\mathbb P^n$.
The following are equivalent:

(i) $\mathcal A$ supports a resonance component of dimension $k-1$ ($k\geq 3$);

(ii) There is a pencil $\mathcal P$ on $\mathbb P^n$ with $k$
{completely} reducible fibers such that $\mathcal A$ is the union of
the zero loci of the linear divisors of all $k$ completely reducible
fibers.

Moreover assume (ii) holds and $F_1,F_2,\ldots,F_k$ are the completely reducible fibers
of $\mathcal P$. Then the 1-forms
$$\omega_i=\frac{dF_i}{F_i}-\frac{dF_k}{F_k},\  i=1,2,\ldots,k-1,$$

form a basis of the corresponding resonance
component of $\mathcal A$.
\end{THM}

Now Theorem \ref{T:5} and our main result Theorem \ref{T:0}
give the upper bounds on the dimensions of resonance components
in terms of ranks of their supports.

\begin{cor}\label{dim}
If  $\mathcal A$ is an arrangement of hyperplanes in $\mathbb P^n$
and $\Sigma \subset \mathbb P( A^1)$ is an  irreducible
component of $\mathcal R^1(\mathcal A)$ then the following
assertions hold:
\begin{enumerate}
\item If $\dim \Sigma > 3$ then the rank of the support of $\Sigma$ is one;
\item If $\dim \Sigma >1$ then the rank  of the support of $\Sigma$ is at most two;
\item If $\dim \Sigma >0$ then the rank  of the support of $\Sigma$  is at most four.
\end{enumerate}
\end{cor}

If the rank of a resonance component is one (for projective arrangement) then
the component is called \emph{local} and is very simple. The result (1) of the corollary
says that the (projective) dimension of non-local component is at most 3.
This in turn implies
that the dimension of $H^1(A,a)$ is at most 3 for $a$ not from a local component.
This result has been proved in \cite{LY2000} for nets and in
\cite{FY} for multinets with all multiplicities of lines equal 1.
Roughly speaking the results (2,3) of Corollary \ref{dim} say that the
non-triviality of resonance varieties is a low-dimensional phenomenon.

Combining Theorem \ref{T:5} with lemma \ref{foliation} we see that
for every resonance component there is a foliation on $\mathbb P^n$.
More directly
the irreducible components of $\mathcal R^1(\mathcal A)$ are
precisely the projectivization of maximal linear  subspaces $E
\subset A^1$  with dimension at least $2$ and isotropic with respect
to the  product $ A^1 \times A^1 \to A^2$, cf. \cite[Corollaries
3.5, 3.7]{LY2000}. In particular, if $E$ is one of these subspaces
then all the homogeneous rational $1$-forms, whose cohomology
classes belong to $E$, are proportional over rational functions
whence correspond to the same foliation. Moreover the  foliation in
question admits a rational first integral $F:\mathbb P^n
\dashrightarrow\mathbb P^1$ of a rather special kind. If we write
$F=\frac{A}{B}$, where $A,B$ are relatively prime homogenous
polynomials, then $sA + tB$ is irreducible for generic $[s:t] \in
\mathbb P^1$
 and the cardinality of the set
\[
\left\{ [s:t] \in \mathbb P^1 \, \vert \, \text{ all irreducible
components of }  s A + t B \text{ have degree one } \right\}
\]
is $\dim E + 1$.

\begin{figure}
\includegraphics[width=4.5cm,height=4.5cm]{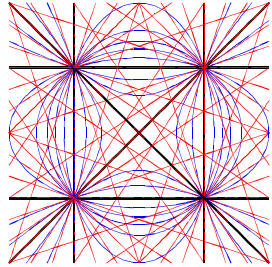}\\
 {\small{\caption{$A_3$ arrangement (in bold) and Bol's $5$-web
  }\label{fbol}}}
\end{figure}

\begin{remark}[Webs associated to Arrangements]\rm
As we have just explained to each arrangement $\mathcal A$ with
$\mathcal R^1(\mathcal A)\not=\emptyset$ we can canonically
associate a finite collection of foliations on $\mathbb P^n$ induced
by the rational maps $F:\mathbb P^n \dashrightarrow \mathbb P^1$,
one for each irreducible component of $\mathcal R^1(\mathcal A)$.
This finite collection forms a global web on $\mathbb P^n$. The
field of web geometry has a venerable history and the present days
are witnessing a lot of activity on webs and their abelian
relations. From this viewpoint one of the key problems, at least
according to  Chern, is the classification of planar webs with the
maximal number of abelian relations that are not algebraizable, cf.
{\cite{chern,bourbaki}}. For a long time the only example that
appeared in the literature was Bol's example. It consists of the
$5$-web formed by four pencils of lines through four generic points
on $\mathbb P^2$ and a pencil of conics through these four points.
It is a rather intriguing fact that the $5$-web, canonically
associated to the Coxeter  arrangement of type A$_3$, is precisely
Bol's web. It corresponds to four resonance components supported on
pencils of lines and one global component. It would be rather
interesting to pursue the determination of the rank of the webs
associated to resonant line arrangements in $\mathbb P^2$.
\end{remark}

\section{Examples and open questions}
\label{S:6}

\subsection{Pencils on $\mathbb P^2$}\label{S:example}

The inequalities of Theorem \ref{T:0} allow pencils on $\mathbb P^2$
 with five completely reducible fibers.
However no such pencils (of full rank) are known.
 This existence problem is not even
settled for the case of nets, cf. \cite[Problem 2]{Ynet}.
The smallest possible example would be a (5,7)-net in $\mathbb P^2$
realizing an orthogonal triple of Latin squares of order 7.
\smallskip

Concerning pencils with four completely reducible fibers
just one example is known -  the Hesse pencil based on the (4,3)-net.
This pencil can be succinctly
described as the pencil generated by a smooth cubic and its Hessian.

\smallskip

Concerning pencils on $\mathbb P^2$ with three completely reducible
fibers,  plenty of them are known, including families with analytic
moduli. The existence of such families with analytic moduli  can be
inferred from the realization result for nets given in \cite[Theorem
4.4]{Ynet}. For {explicit examples} one can consider hyperplane
sections of the examples presented in \ref{S:P3}.

\subsection{Pencils  on $\mathbb P^3$}\label{S:P3}

\begin{figure}
\includegraphics[width=4.5cm,height=4.5cm]{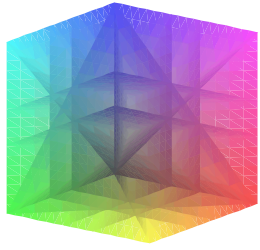}\\
  {\small{\caption{Arrangement of $12$ hyperplanes on $\mathbb P^3$  associated to
 $D_4$
  }\label{p4}}}
\end{figure}

Let $d$ be a positive integer and consider the pencil on $\mathbb
P^3$ generated by
\begin{eqnarray*}
F_d &=&  \left( x_0 ^d - x_1^d  \right)  \left( x_2^d - x_3^d
\right)
 \\
G_d &=&  \left( x_0 ^d - x_2^d  \right)  \left( x_1^d - x_3^d
\right) \,.
\end{eqnarray*}
Then
\[
 F_d -  G_d = \left( x_0 ^d - x_3^d  \right)  \left( x_2^d - x_1^d
 \right) \,.
\]
The arrangements $\mathcal A_d$ corresponding to these pencils have many interesting
properties.
Each $\mathcal A_d$ consists of $6d$ hyperplanes that are the reflecting hyperplanes
for the monomial group $G(d,d,4)$ generated by complex reflections (e.g., see \cite{OT}).
For $d=1,2$ the group is the Coxeter group of type $A_3$ or $D_4$ respectively.
For all $d\ge 2$ the arrangements $\mathcal
A_d$ are essential and carrying each a global irreducible resonance component of dimension one.

Moreover these arrangements are (3,d)-nets
in $\mathbb P^3$ (for $d\ge 2$).
This implies that the intersection of $\mathcal A_d$
with a generic plane gives a net in $\mathbb P^2$. On the other hand, the intersection
with a special plane (say, the one defined by $x_3=0$) gives the family of multinets mentioned at the end of Section 3 whence representing these multinets as limits of nets.

Let us look closer at the combinatorics of $\mathcal A_d$, i.e., at the corresponding
Latin square (see \cite{Ynet}). For each $\zeta$ such that $\zeta^d=1$
denote by $H_{i,j}(\zeta)$ the hyperplane defined by $x_i=\zeta x_j$ ($1\leq i<j\leq 4$). Then identify
the collection of hyperplanes corresponding to the linear divisors of $F_d$, $G_d$. and $F_d-G_d$
via
$$a_{\zeta}=H_{1,2}(\zeta^{-1})=H_{1,3}(\zeta)=H_{2,3}(\zeta),$$
and
$$b_{\zeta}=H_{3,4}(\zeta)=H_{2,4}(\zeta)=H_{1,4}(\zeta).$$
After the identification the Latin square corresponding to
$\mathcal A_d$ is the multiplication table of the dihedral group
$D_d$ with $a_{\zeta}$ forming the cyclic subgroup of order $d$ and
$\{b_{\zeta}\}$ being the complementary set of involutions.
Intersecting $\mathcal A_d$ for $d\geq 3$ with a general plane we
obtain a series of 3-nets in $\mathbb P^2$ realizing noncommutative
groups (cf. \cite{Ynet}). In particular these nets are not
algebraizable. We remark that another non-algebraizable example of a
3-net in $\mathbb P^2$ has been found by J. Stipin in \cite{St}. He
has exhibited a (3,5)-net that does not realize $\mathbb Z_5$.

Finally (as we learned from \cite{St}) the general fibers of the
above pencil for $d=2$ were studied by R. M. Mathews in \cite{Ma}
under the name of \emph{desmic} surfaces. Thus for arbitrary $d$
these fibers can be considered as generalizations of desmic
surfaces.

\smallskip

The foliations  $\mathcal F_d$ induced by $\omega_d= F_d dG_d - G_d
dF_d$ have degree $4d-2$. Thus the bound of Proposition \ref{P:para}
is attained by $\mathcal F_2$ and no other foliation in the family.

\subsection{Pencils  on $\mathbb P^4$}\label{S:P4}
We do not know any example of a pencil on $\mathbb P^4$ with three completely
reducible fibers that is not a linear
pullback of a pencil from a smaller dimension. One can deduce from careful reading of
the proof of Proposition \ref{red} that the degree of such a pencil must be at least 10.

\section{Pencils and characteristic varieties}
\label{S:7}

Let $\mathcal A$ be an arrangement in $\mathbb P^n$ and $M$ its
{complement}.
 If for every $\rho \in
\mathrm{Hom}\left(\pi_1( M ), \mathbb C^*\right)$ we write $\mathcal
L_{\rho}$ for the associated rank one local system then  the
characteristic varieties $\mathcal V^l(M)$ are defined as follows
\[
\mathcal V^l(M) = \left \{ \rho \in \mathrm{Hom}(\pi_1( M , \mathbb
C^*) \, \vert \, \mathrm H^l(M,\mathcal L_{\rho}) \neq 0 \right \}.
\]
 If $\mathcal A'\subset \mathcal A$ is a
subarrangement with complement $M'$ then the inclusion
of $M$ into $M'$ induces an inclusion of
$\mathrm{Hom}(\pi_1( M' ), \mathbb
C^*)$ into $\mathrm{Hom}(\pi_1( M) , \mathbb
C^*)$ and also of $\mathcal V^1(M')$ into $\mathcal V^1(M)$.
As in the case of resonance varieties, the support of an irreducible component $\Sigma$
of $\mathcal V^1(M)$ is the smallest subarrangement $\mathcal A'$ such that
$\Sigma \subset \mathcal V^1(\mathcal A')$.

The main result of \cite{CS1999}  implies that the projectivization
of the tangent cone of $\mathcal V^l(M)$ at the trivial
representation is isomorphic to $\mathcal R^l(\mathcal A)$.
Thus if $\Sigma \subset \mathcal V^1(M)$ is an irreducible component of
dimension $d>0$ containing $1 \in \mathrm{Hom}(\pi_1( M) , \mathbb
C^*)$ then the
projectivization of its tangent space is an irreducible component of
the resonance variety.

There exist translated components of characteristic varieties, i.e., those that
do not contain the trivial representation, see \cite{Su}.
In a recent preprint  Dimca  \cite{Dimca} found more precise properties of
shifted components and clarified  the link
between the positive dimensional irreducible components of $\mathcal
V^1(M)$ and pencils of hypersurfaces.
We invite the reader to consult \cite{Dimca} for a more extensive description. Here we will
recall just what is strictly necessary for our purposes.

If $\Sigma$ is a translated component of dimension
at least two then after translating it to $1$ one
obtains an irreducible component of $\mathcal V^1(M)$ through
$1$.

If $\Sigma$ is a translated component of dimension one then there
exists a pencil of hypersurfaces with generic irreducible fiber and
exactly two completely reducible fibers with support contained in
the support of $\Sigma$. The support of $\Sigma$ is the union of the
hyperplanes that appear as components of {fibers} of the pencil.
Moreover there is at least one extra fiber such that its components
are either hyperplanes in the support of $\Sigma$ or non-reduced
hypersurfaces. Combining these results with our methods we can prove
the upper bound for the dimension of characteristic components.

The proof of the following theorem repeats almost verbatim the proof
of Theorem \ref{T:1} using now Proposition \ref{L:like} in full. So
we omit {its} proof.

\begin{THM}\label{T:6}
Let $\mathcal P$ be an essential pencil of hypersurfaces on $\mathbb P^n$
with irreducible generic fiber, two completely reducible fibers, and
a third fiber that is a product of linear forms and non-reduced polynomials.
Then the Gauss map of the associated
foliation $\mathcal F$ is non-degenerate.
\end{THM}

Now using again Dimca's results and Theorem \ref{T:6}
we can prove an analogue of Corollary \ref{dim} for
characteristic varieties.

\begin{THM}\label{T:7}
If  $\mathcal A$ is an arrangement of hyperplanes in $\mathbb P^n$
and $\Sigma \subset \mathcal V^1(M)$ is an  irreducible
component then the following
assertions hold:
\begin{enumerate}
\item If $\dim \Sigma >4$ then the rank of the support of $\Sigma$ is one;
\item If $\dim \Sigma >2$ then the rank  of the support of $\Sigma$ is at most two;
\item If $\dim \Sigma >1$ then the rank  of the support of $\Sigma$  is at most four;
\item If $\dim \Sigma >0$ then the rank  of the support of $\Sigma$  is at most six.
\end{enumerate}
\end{THM}

(Notice that the difference with Corollary \ref{dim}
in the dimensions of $\Sigma$ is due to the projectivization in the corollary.)

\begin{proof}
The statements (1)-(3) follow immediately from Corollary \ref{dim}
and the relations between resonance and characteristic components.
In order to prove (4) it suffices to
show that $n<7$ if  there exists a pencil $\mathcal P$ of
hypersurfaces on $\mathbb P^n$ with completely reducible generators
$F$ and $G$ inducing a full rank arrangement and at least one extra
fiber, say $K=F-G$, that can be written as
\[
  K=\widetilde U \cdot \widetilde V
\]
where $\widetilde U\in \mathbb C[x_0, \ldots,x_n]$ is  a product of
linear forms and $\widetilde V\in \mathbb C[x_0, \ldots,x_n]$ is a
product of non-trivial powers of irreducible polynomials of degree
at least two. We  denote by $U,V$ the reduced polynomial with the
same zero set of $\widetilde U, \widetilde V$ and point out that
\begin{equation}\label{E:trivial}
2\deg \left( \frac{\widetilde V}{V} \right) \ge \deg (\bar V) \, .
\end{equation}

If $\widetilde Q \in \mathbb C[x_0, \ldots,x_n]$ denotes the product
$FG$ and $Q$ denotes the reduced polynomial with the same zero set
then
\[
 \omega = \frac{U}{\widetilde U}\frac{V}{\bar V}\frac{Q}{\widetilde Q} \omega_0
\]
is a homogeneous polynomial $1$-form defining the foliation $\mathcal F$ associated to
$\mathcal P$. In
particular
\[
  \deg(\mathcal F) \le \deg(\widetilde Q) - 2 - \deg\left(\frac{\widetilde U \widetilde V \widetilde Q}{UVQ}\right) .
\]

  Theorem \ref{T:7} implies
that the Gauss map of $\mathcal F$ is non-degenerate. The property
(iv) of the forms proved in section \ref{S:2} can be immediately
generalized to $\omega$ and it implies that  $(\widetilde U
\widetilde Q)^{n-1}$ divides
\[
  \left(\frac{\widetilde U \widetilde Q}{U Q}\right)^{n-1}\det \left(\frac{\partial a_i}{\partial x_j}
  \right) \, ,
\]
where $\omega = \sum a_i dx_i$. We obtain that
\begin{eqnarray*}
(n-1) \deg(\widetilde U \widetilde Q)) &\le& (n+1) \left( \deg(
\widetilde Q) -2 -\deg\left(\frac{\widetilde U \widetilde V
\widetilde Q}{UVQ}\right)\right) +
(n-1)\deg\left(\frac{\widetilde U \widetilde Q}{UQ}\right) \\
&=& (n+1) \left(\deg( \widetilde Q) -2 \right)  -
2\deg\left(\frac{\widetilde U \widetilde Q}{U Q}\right) -(n+1)
\deg\left(\frac{\widetilde V}{V}\right) .
\end{eqnarray*}
Therefore, if we suppose that $n\ge 3$,  delete the term $-
2\deg\left(\frac{\widetilde U \widetilde Q}{U Q}\right)$ and use
(\ref{E:trivial}) then
\begin{eqnarray*}
(n-1)  \deg(\widetilde Q)  &<&  (n+1) \left( \deg(\widetilde Q) -2
\right) -\frac{n+1}{2}
\deg\left(\widetilde V\right) -(n-1) \deg ( \widetilde U)  \\
&\le& (n+1) \left( \deg(\widetilde Q) -2 \right) -\frac{n+1}{2}
\deg\left(\widetilde U\widetilde V\right) < \frac{3(n+1)}{4}
\deg\left(\bar Q\right).
\end{eqnarray*}
In particular $n < 7$ and the statement follows.
\end{proof}

We do not know if the bounds given in the theorem are sharp.

We also do not know if there are restrictions on the rank of the
support of the zero-dimensional characteristic varieties.

\end{document}